\documentclass{amsart}[12pt]

\theoremstyle{definition}

\theoremstyle{remark}

\newtheorem{Definition}{\bf Definition}[section]
\newtheorem{Thm}[Definition]{\bf Theorem}

\newtheorem{Cor}[Definition]{\bf Corollary}
\newtheorem{Prop}[Definition]{\bf Proposition}
\newtheorem{Example}[Definition]{\bf Example}

\numberwithin{equation}{section}

\reversemarginpar

\newcommand{\nn}{\nonumber}
\newcommand{\no}{\noindent}
\newcommand{\realpart}{\mathop{\rm Re}\nolimits}

\newcommand{\ba}{\begin{eqnarray}}
\newcommand{\ea}{\end{eqnarray}}

\newcommand{\ione}{\int_{0}^{1}}

\newcommand{\allC}{\mathbb{C}}
\newcommand{\allN}{\mathbb{N}}
\newcommand{\nnN}{\mathbb{N}_{0}}

\newcommand{\stvictor}{}

\newcommand{\barzeta}{{\bar\zeta}}

\begin{document}

\title[The evaluation of Tornheim double sums. Part 1] {The evaluation of
Tornheim double sums. Part 1}

\author{Olivier Espinosa}
\address{Departamento de F{\'\i}sica,
Universidad T{{\'e}}c. Federico Santa Mar{\'\i}a, Valpara{\'\i}so, Chile}
\email{olivier.espinosa@usm.cl}

\author{Victor H. Moll}
\address{Department of Mathematics,
Tulane University, New Orleans, LA 70118}
\email{vhm@math.tulane.edu}

\subjclass{Primary 33}

\date{\today}


\begin{abstract}
We provide an explicit formula for the Tornheim double series in terms of
integrals involving the Hurwitz zeta function. We also study the limit when
the parameters of the Tornheim sum become natural numbers, and show that
in that case it can be expressed in terms of definite integrals of triple
products of Bernoulli polynomials and the Bernoulli function
$A_k (q): = k\zeta '(1 - k,q)$.

\end{abstract}

\maketitle

\section{Introduction} \label{S:intro}

The function
\ba
T(a,b,c) & = & \sum_{n=1}^{\infty} \sum_{m=1}^{\infty} \frac{1}{n^{a}
\, m^{b} \, (n+m)^{c} }, \quad a, \, b, \, c \in \mathbb{C},
\label{seriestz}
\ea
\no
was introduced by Tornheim in \cite{torn}.  We
provide here an analytic expression for $T(a,b,c)$ in terms of
the integrals
\begin{align}
\label{intI}
I(a,b,c) & = \ione \zeta(1-a,q) \zeta(1-b,q) \zeta(1-c,q) \, dq\\
\intertext{and}
\label{intJ}
J(a,b,c) & = \ione \zeta(1-a,q) \zeta(1-b,q) \zeta(1-c,1-q) \, dq.
\end{align}
\no
Here $\zeta(z,q)$ is the Hurwitz zeta function,
\ba
\zeta(z,q) & = & \sum_{n=0}^{\infty} \frac{1}{(n+q)^{z}},
\label{hurzetadef}
\ea
\no
defined for  $z \in \mathbb{C}$ and
$q \neq 0, \, -1, \, -2, \cdots $. The series (\ref{hurzetadef}) converges for
$\realpart{z} > 1$ and $\zeta(z,q)$ admits a meromorphic extension  to the
complex plane with a single pole at $z=1$ as its only singularity.

In the case where the parameters $a,b,c$ in \eqref{seriestz} are
positive integers, the Tornheim sum can be expressed in terms of
the Riemann zeta function
\ba
\zeta(z) = \zeta(z,1) = \sum_{n=1}^{\infty} \frac{1}{n^{z}},  \label{riemannzeta}
\ea
its derivatives, and integrals related to the families \eqref{intI} and
\eqref{intJ}, as given in Theorem \ref{allsums} below.

We use the notation
\ba
\barzeta(z,q):=\zeta(1-z,q),
\ea
defined for $z\neq 0$ and $q \neq 0, \, -1, \, -2, \cdots$.

\medskip

The results presented here  are a
continuation of \cite{esmo1,esmo2} where we
have provided many explicit evaluations of definite integrals containing
$\zeta(z,q)$ in the integrand. For instance, if $\realpart{a} > 0$, then
\ba
\int_{0}^{1} \barzeta(a,q) \, dq & = &  0,
\label{intzeta1}
\ea
\no
and, for $\realpart{a} > 1, \, \realpart{b} > 1$, we have
\begin{align}
\label{twozetaa}
\int_{0}^{1} \barzeta(a,q) \barzeta(b,q) \, dq & =
\frac{2 \Gamma(a) \Gamma(b)}{(2 \pi)^{a+b}}
\, \zeta(a+b) \cos\left( \frac{\pi}{2}(a-b) \right)
\stvictor
\\
\intertext{and}
\label{twozeta1a}
\int_{0}^{1} \barzeta(a,q) \barzeta(b,1-q) \, dq & =
\frac{2 \Gamma(a) \Gamma(b)}{(2 \pi)^{a+b}}
\, \zeta(a+b) \cos\left( \frac{\pi}{2}(a+b) \right).
\stvictor
\end{align}
\no
Lerch's evaluation  \cite{ww},
\ba
\frac{d}{dz} \zeta(z,q) \Big{|}_{z=0} & = & \ln \Gamma(q) - \ln \sqrt{2 \pi},
\stvictor
\ea
\no
yields integrals involving the
loggamma function. For instance,
\ba
L_{1}= \int_{0}^{1} \ln \Gamma(q) \, dq & = & \ln \sqrt{2 \pi}
\label{l1}
\stvictor
\ea
\no
and
\ba
L_{2}= \int_{0}^{1} \ln^{2} \Gamma(q) \, dq & = & \frac{\gamma^{2}}{12} +
\frac{\pi^{2}}{48} + \frac{\gamma \, L_{1}}{3} + \frac{4}{3}L_{1}^{2}
- \frac{A \, \zeta'(2)}{\pi^{2}} + \frac{\zeta''(2)}{2 \pi^{2}}
\label{l2}
\stvictor
\ea
\no
with
\ba
A & = &  \gamma + \ln 2 \pi.
\label{A}
\stvictor
\ea

We expect that the methods developed here will provide analytic
expressions for the constant
\ba
L_{3} & = & \int_{0}^{1} \ln^{3} \Gamma(q) \, dq.
\ea

\medskip

The series $T(a,b,c)$, for $a, \, b, \, c \in \mathbb{R} - \mathbb{N}$, is
given in Theorem \ref{thm36} in terms of integrals (\ref{intI}) and
(\ref{intJ}).  The evaluation of the Tornheim series for integer values of
the parameters are expressed in terms of some definite integrals:

\begin{Thm}
\label{allsums}
The Tornheim sums $T(n_{1},n_{2},n_{3})$ can be expressed as a finite expression
of the Riemann zeta function, its derivatives and the integrals
\ba
K_{m,n} & = & \ione \psi^{(-m)}(q) B_{n}(q) \ln \Gamma(q)\, dq, \nn \\
K^{*}_{m,n} & = & \ione  \psi^{(-m)}(1-q) B_{n}(q)\ln \Gamma(q) \, dq,
\nn \\
Z_{m,n} & = &
\ione \psi^{(-m)}(q)  \psi^{(-n)}(q) \ln \Gamma(q) \, dq, \nn \\
Z_{m,n}^{*} & = &
\ione \psi^{(-m)}(q)  \psi^{(-n)}(1-q) \ln \Gamma(q) \, dq. \nn
\ea
\end{Thm}

Here $B_{n}(q)$ is the Bernoulli
polynomial, $\Gamma(q)$ is the classical gamma function and
$\psi^{(-m)}(q) = \psi(-m,q)$, where
\ba
\psi(z,q) & = & e^{- \gamma z} \frac{\partial}{\partial z}
\left[ e^{\gamma z} \frac{\zeta(z+1,q)}{\Gamma(-z)} \right]
\label{psi-def1}
\ea
\no
is the generalization of the polygamma function introduced
by the authors in  \cite{esmo4}, with $z$ an arbitrary complex variable.
Some properties of $\psi(z,q)$ are given in Appendix 2.
The closed form evaluation of the integrals in Theorem \ref{allsums}
will be discussed in a future paper.
The proof of Theorem \ref{allsums} is given in Sections \ref{S:limit}
and \ref{S:Q-integrals}.

\medskip

The series (\ref{seriestz}) converges for
$\realpart{a}, \, \realpart{b}, \, \realpart{c} > 1$.
Matsumoto \cite{matsu1} showed that it can be continued as
a meromorphic function to $\mathbb{C}^{3}$, with all its singularities
located on the subsets of $\mathbb{C}^{3}$ defined by one of the equations
\ba
\nn
a+c=1 - l, \, b+c = 1- l, \, a+b+c=2 \quad \text{ with } \, l \in \nnN.
\ea

The literature contains many techniques to evaluate some particular
cases of $T(a,b,c)$. For instance, the case $c= 0$ is evaluated simply as
\ba
T(a,b,0) & = &  \zeta(a) \, \zeta(b).  \nn
\stvictor
\ea
\no
The elementary identity
\ba
T(a,b-1,c+1) + T(a-1,b,c+1) & = & T(a,b,c)
\label{elem00}
\stvictor
\ea
\no
and the symmetry rule
\ba
T(a,b,c) &  = & T(b,a,c) \label{symmetryrule}
\stvictor
\ea
\no
has been used by Huard et al \cite{huard} to give
the explicit expression
\begin{equation*}
T(a,b,c) =  \sum_{i=1}^{a} \binom{a+b-i-1}{a-i} T(i,0,N-i)  +
\sum_{i=1}^{b} \binom{a+b-i-1}{b-i} T(i,0,N-i)
\stvictor
\end{equation*}
\no
in the case that both $a$ and $b$ are positive integers.
Here $N = a+b+c$. If $N$ is an odd positive integer greater than $1$,
then the sum $T(i,0,N-i)$ is evaluated as
\begin{align*}
T(i,0,N-i) = (-1)^{i} \sum_{j=0}^{\lfloor{ (N-i-1)/2 \rfloor} }
&\binom{N-2j-1}{i-1} \zeta(2j) \zeta(N-2j) \\
+ (-1)^{i} \sum_{j=0}^{\lfloor{ i/2 \rfloor} }
&\binom{N-2j-1}{N-i-1} \zeta(2j) \zeta(N-2j) +  \zeta(0) \zeta(N)
\stvictor
\end{align*}
\no
The evaluation of $T(i,0,N-i)$ in the case
$N$ even remains open.  The techniques introduced in this paper have allowed us
to evaluate the sum $T(a,0,c)$ in terms
of integrals similar to the ones discussed here.
Details will appear in  \cite{esmo11}.

The multiple zeta value, also called Euler sums, are defined by
\ba
\zeta(i_{1},i_{2}, \ldots, i_{k} ) & = & \sum \frac{1}{n_{1}^{i_{1}}
\, n_{2}^{i_{2}} \cdots n_{k}^{i_{k}} }
\ea
\no
where the sum extends over $n_{1} > n_{2} > \cdots > n_{k} > 0$. The
sum $T(a,0;c)$ is $\zeta(c,a)$. A general introduction to these sums is
provided in Chapter 3 of Borwein-Bailey-Girgenson \cite{borbaigir}.

\medskip

The identities of Tornheim \cite{torn}  for
$T(n_{1},n_{2},n_{3})$ are based on an
elementary identity for series: let $f$ be monotone decreasing and $f(x) \to c$ as
$x \to \infty$ and define
\ba\nn
\varphi(n,m;f) & = & \frac{f(m)}{n(n+m)} + \frac{f(n)}{m(n+m)} -
\frac{f(n+m)}{nm}.
\ea
\no
Then
\ba\nn
\sum_{n,m=1}^{\infty} \varphi(n,m;f)  & = & 2 \sum_{r=1}^{\infty}
\frac{f(r) - c}{r^{2}}.
\ea
\no
The special case $f(x) = 1/x$ yields
\ba
T(1,1,1) & = & 2 \zeta(3)
\stvictor
\ea
\no
and $f(x) = 1/x^{a-2}$ produces the relation
\ba\nn
2T(a-2,1,1) - T(1,1,a-2) & = & 2 \zeta(a).
\stvictor
\ea
\no
Among the many evaluation presented in \cite{torn} we mention
\begin{align*}
T(1,1,a-2) &= (a-1) \zeta(a) - \sum_{i=2}^{a-2} \zeta(i) \zeta(a-i),\stvictor\\
T(a-2,1,1) &= \frac{1}{2}T(1,1,a-2) + \zeta(a)\\
\intertext{and}
T(1,0,a-1) &= \frac{1}{2} T(1,1,a-2).
\end{align*}

\no
Subbarao and Sitaramachandrarao \cite{sita} give
\begin{align}
\label{T-sym-even}
T(2n,2n,2n) & =\frac{4}{3} \sum_{i=0}^{n} \binom{4n-2i-1}{2n-1}
\zeta(2i) \zeta(6n-2i),\stvictor
\\
\intertext{and this is complemented by Huard et al. \cite{huard} with}
\label{T-sym-odd}
T(2n+1,2n+1,2n+1) & = -4 \sum_{i=0}^{n} \binom{4n-2i+1}{2n} \zeta(2i)
\zeta(6n-2i+3).\stvictor
\end{align}

\medskip

Boyadzhiev \cite{boya1,boya2} has given elementary proofs of an expression
for $T(a,b,c)$ in terms of the function
\ba\nn
S(r,p) & = & \sum_{n=1}^{\infty} \frac{H_{n}^{(r)}}{n^{p}}.
\ea
\no
Here $H_{n}^{(r)} = 1^{-r} + 2^{-r} + \ldots + n^{-r}$ is the generalized
harmonic number. In \cite{boya2} the author establishes recurrences for the
sums $S(r,p)$ that permit to express them as products of zeta values in the
case $N = r+p$ odd.

\medskip

Tornheim double sums and other related ones appear as special
cases of the zeta function $\zeta_{\mathfrak{g}}(s)$ of a semi-simple
Lie algebra defined as
\ba\nn
\zeta_{\mathfrak{g}}(s) & = & \sum_{\rho} \text{dim}(\rho)^{-s}
\ea
\no
where the sum is over all the finite-dimensional representations
of $\mathfrak{g}$.
Zagier \cite{zagier} states that the special case
$\mathfrak{g} = \mathfrak{s} \mathfrak{l}(3)$
yields (\ref{T-sym-even}).  The nomenclature
for $T(a,b,c)$ is not standard: Zagier \cite{zagier} and also
Crandall and Buhler \cite{crandall} call $T(a,b,c)$ the {\em Witten zeta
function}. Tsumura \cite{tsumura1} has evaluated some special cases of the
sum
\ba\nn
W(p,q,r,s) & = & \sum_{m, n =1}^{\infty}
\frac{1}{m^{p} n^{q} (m+n)^{r} (m+2n)^{s} }
\ea
\no
under the parity restriction $p+q+r+s$ is odd.  This is the Witten sum
corresponding to $SO(5)$. \\

These sums also have appeared in connection with knots and
Feynman diagrams, see \cite{kreimer} for details. \\

The Bernoulli function
\ba
A_{k}(q) & = & k \zeta'(1-k,q), \quad k \in \mathbb{N}
\ea
\no
introduced in \cite{esmo4}, plays an important role in the evaluations
presented here. Adamchik \cite{adambarnes} proved the identity
\ba
\zeta'(1-k,q) & = & \zeta'(1-k) +
\sum_{j=0}^{k-1} (-1)^{k-1-j} j! Q_{j,k-1}(q) \, \ln \Gamma_{j+1}(q)
\ea
\no
where
\ba
Q_{k,n}(q) & = & \sum_{j=k}^{n}  (1-q)^{n-j} \binom{n}{j}
\left\{ \begin{matrix}  j  \\  k \end{matrix}  \right\}
\ea
\no
is the Stirling polynomial and the generalized gamma function $\Gamma_{n}(q)$
is defined inductively via
\ba
\Gamma_{n+1}(q+1) & = & \frac{\Gamma_{n+1}(q) }{\Gamma_{n}(q)} \label{gammagen} \\
\Gamma_{1}(q) & = & \Gamma(q) \nn \\
\Gamma_{n}(1) & = & 1. \nn
\ea

\newpage

\no{\bf Notation}:

\medskip

{\allowdisplaybreaks
\begin{tabular}{cl}
$\Gamma(q)$ & is the gamma function, \\
$\zeta(z)$ & is the Riemann zeta function defined in \eqref{riemannzeta}, \\
$\zeta(z,q)$ & is the Hurwitz zeta function defined   in \eqref{hurzetadef},\\
$\barzeta(z,q)$ & is a shorthand for $\zeta(1-z,q)$, \\
$\zeta_{\pm}(z,q)$ & denotes the combination $\zeta(z,q) \pm \zeta(z,1-q)$, \\
$\psi(z,q)$ & is the generalized polygamma function, defined in \eqref{psi-def1}, \\
$\psi^{(-n)}(q)$ & is the balanced negapolygamma function, defined in \eqref{bal-negapolygamma}, \\
$A_{n}(q)$ & is the Bernoulli function, defined in \eqref{def-A}, \\
$B_{n}(q)$ & is the Bernoulli polynomial of degree $n$, defined in \eqref{bergen1} and \eqref{bergen3}, \\
$B_{n}$ & is the $n$-th Bernoulli number,\\
$H_{n}$ & is the $n$-th harmonic number,  $H_n = 1 + \tfrac{1}{2} + \cdots + \frac{1}{n}$, \\
$h_{n}$ & is a shorthand for $H_{n-1}$, \\
$\gamma$ & is Euler's constant, \\
$A$ & equals $\gamma + \ln 2 \pi$, \\
$A_{\pm}$  & equals $ A^{2} \pm \frac{\pi^{2}}{4}$.
\end{tabular}
}

\medskip

\section{The main identity} \label{S:identity}

We now provide an analytic expression for the Tornheim
double series  $T(a,b,c)$ in terms of the integrals
(\ref{intI}) and (\ref{intJ}). The
analysis of its behavior as the
parameters become integers is described in Section \ref{S:limit}.  The
proof employs the Fourier representation for $\barzeta(z,q)$:
\begin{multline}
\barzeta(z,q) = \frac{2 \Gamma(z)}{(2 \pi)^{z}}
\left[ \cos \left( \frac{\pi z}{2} \right) \sum_{n=1}^{\infty}
\frac{\cos(2 \pi q  n)}{n^{z}} +
\sin \left( \frac{\pi z}{2} \right) \sum_{n=1}^{\infty}
\frac{\sin(2 \pi q  n)}{n^{z}} \right],
\label{fourier0}
\end{multline}
\no
valid for $\realpart{z} > 1$ and $0 < q \le 1$, given in \cite{berndt}.

From (\ref{fourier0}) we obtain
\ba
2 f_{c}(z)  \sum_{n=1}^{\infty} \frac{ \cos 2 \pi q n}{n^{z}} & = &
\barzeta(z,q) + \barzeta(z,1-q) \label{zetacos}
\ea
\no
and
\ba
2 f_{s}(z) \sum_{n=1}^{\infty} \frac{ \sin 2 \pi q n}{n^{z}} & = &
\barzeta(z,q) - \barzeta(z,1-q), \label{zetasin}
\ea
\no
where
\begin{equation}
f_{c}(z) = \frac{2 \Gamma(z)}{(2 \pi)^{z}} \cos \left( \frac{\pi z}{2} \right)
\quad \text{ and }\quad
f_{s}(z) = \frac{2 \Gamma(z)}{(2 \pi)^{z}} \sin \left( \frac{\pi z}{2} \right).
\label{deffcfs}
\end{equation}

\medskip

For a function $h(a,b,c)$ we denote
\begin{align}
h^{sym}(a,b,c) & = h(a,b,c) + h(b,c,a) + h(c,a,b) \label{hplus}\\
\intertext{and}
h^{nsym}(a,b,c) & = -h(a,b,c) + h(b,c,a) + h(c,a,b). \label{hminus}
\end{align}

\medskip

\begin{Prop}
\label{thm31}
Let $a, \, b, \, c \in \mathbb{R}$. Then $T(a,b,c)$ satisfies the relations
\begin{align}
\label{T-sym}
f_{c}(a)f_{c}(b)f_{c}(c) T^{sym}(a,b,c)  & =
I(a,b,c) + J^{sym}(a,b,c)\\
\intertext{and}
\label{T-nsym}
f_{s}(a)f_{s}(b)f_{c}(c) T^{nsym}(a,b,c)  & =
I(a,b,c) - J^{nsym}(a,b,c),
\end{align}
\no
where $f_{c}$ and $f_{s}$ are defined in (\ref{deffcfs}). \\
\end{Prop}
\begin{proof}
Multiply three series of cosine type in (\ref{zetacos}) to obtain that
\ba
8f_{c}(a) f_{c}(b) f_{c}(c)
\sum_{n=1}^{\infty} \sum_{m=1}^{\infty} \sum_{r=1}^{\infty}
\frac{1}{n^{a} \, m^{b} \, r^{c} }
\int_{0}^{1} \cos(2 \pi q n) \, \cos(2 \pi q m) \, \cos(2 \pi q r) \, dq
\nn
\ea
equals
\begin{multline}
\nn
\int_{0}^{1} \left[\barzeta(a,q) + \barzeta(a,1-q)\right]
\left[\barzeta(b,q) + \barzeta(b,1-q)\right]
\left[\barzeta(c,q) + \barzeta(c,1-q)\right] \, dq.
\end{multline}
\no
The identities
\ba
4 \cos(u) \cos(v) \cos(w)  & = &
\cos(u+v+w) + \cos(u+v-w)  \nn \\
 & & \; + \, \cos(u-v+w)+ \cos(u-v-w)
\nn
\ea
\no
and
\ba\nn
\int_{0}^{1} \cos(2 \pi q j) \, dq & = & \begin{cases}
                                          0 \text{  if } j \neq 0 \\
                                          1 \text{  if } j  = 0,
                                         \end{cases}
\ea
\no
reduce the left-hand side to $ 2f_{c}(a) f_{c}(b) f_{c}(c) T^{sym}(a,b,c)$.
To complete the proof of the first identity, we expand the products of zeta
functions to write the integral as a sum of eight different integrals, which can
be reduced to the right-hand side of \eqref{T-sym}, by selectively performing the
the change of variable $q \to 1-q$ in half of them.

The second identity is obtained by considering the only other nonvanishing
triple product integral, namely, that of
$\sin(2 \pi q n) \sin(2 \pi q m) \cos(2 \pi q r)$.
\end{proof}

\medskip

The case of Proposition \ref{thm31} in which the parameters are
integers will be our main interest in this paper. When the argument $z$ of
the function $\barzeta(z,q)$ is a positive integer $n$, this function reduces
to a Bernoulli polynomial,
\ba\label{Hurwitz-to-Bernoulli}
\barzeta(n,q)=-\frac{1}{n}B_n(q).
\ea
In this case, due to the reflection property of the Bernoulli polynomials,
\ba\label{Bernoulli-reflection}
B_k (1 - q) = ( - 1)^k B_k (q),
\ea
the function $\barzeta(n,1-q)$ reduces simply to $\barzeta(n,q)$ up to a sign:
\ba
\barzeta(n,1-q)=(-1)^n \barzeta(n,q),
\ea
so that the $J$-type integrals reduce to $I$-type integrals.\\

Unfortunately, since the functions $f_c(n)$ and $f_s(n)$ vanish for $n$ odd and
$n$ even respectively, the identities of Proposition \ref{thm31} for integer parameters
$(a,b,c)=(n_1,n_2,n_3)$ are trivial except only in two cases:\\
\begin{enumerate}
    \item $n_1,n_2,n_3$ are all even,
    \item $n_1,n_2$ are both odd, and $n_3$ is even.
\end{enumerate}
The first case is of special interest. It appears in \cite{sita} as
{\em the reciprocity relation for a class of Tornheim series}.

\begin{Cor}
\label{cor33}
Let $n_{1}, \, n_{2}, \, n_{3} \in \mathbb{N}$ be even. Then
\begin{multline}
T^{sym}(n_{1},n_{2},n_{3}) =
(-1)^{(n_{1}+n_{2}+n_{3})/2}\frac{(2 \pi)^{n_{1}+n_{2}+n_{3}}}
{2(n_{1}-1)! (n_{2}-1)! (n_{3}-1)!} I(n_{1},n_{2},n_{3}).
\end{multline}
\end{Cor}

\medskip

\begin{Cor}
\label{cor65}
Let $n \in \mathbb{N}$. Then
\begin{align}
\label{T-sym-even-2}
T(2n,2n,2n) =  \frac{1}{3}(-1)^{n} \, (2 \pi)^{6n}
\sum_{k=0}^{n} \binom{4n-2k-1}{2n-1}
\frac{B_{2k} B_{6n-2k}}{(2k)!(6n-2k)!}.
\end{align}
\end{Cor}
\begin{proof}
Corollary \ref{cor33} and \eqref{Hurwitz-to-Bernoulli} yield
\ba\nn
T(2n,2n,2n) =
(-1)^{n+1}\frac{(2 \pi)^{6n}}
{6(2n)!^3} \ione B_{2n}(q)^{3} \, dq .
\ea
\no
The value of the integral is given by Carlitz \cite{carlitz}.
It can also be obtained directly from the formula for $B_n(q)^3$ given in
Appendix 1.
\end{proof}
\no
Formula \eqref{T-sym-even-2} agrees with formula \eqref{T-sym-even} on account
of the relation
\ba
\label{rel-zeta-bernoulli-even-n}
\zeta(2k) = \frac{(-1)^{k+1} \, (2 \pi)^{2k} B_{2k} }{2 (2k)!},
\ea
valid for $k \in \nnN$.

\medskip

We now present an analytic expression for the Tornheim double series, valid
for non-integer values of the parameters. \\

\begin{Thm}
\label{thm36}
Let $a,b,c \in \mathbb{R}$ and define
\ba
\lambda(z) & = & \frac{\Gamma(1-z)}{(2 \pi)^{1-z}} =
\frac{\pi}{(2 \pi)^{1-z} \, \Gamma(z) \, \sin \pi z}.
\label{lambdadef}
\ea
\no
For $a,b,c \not \in \mathbb{N}$ we have
\begin{multline}\label{T-explicit}
T(a,b,c) =  4 \lambda(a) \lambda(b) \lambda(c) \sin(\pi c/2)
\Big[ \cos\left(\tfrac{\pi}{2}(a-b)\right) \left[ J(c,a,b)+J(c,b,a) \right]
\\
 -  \cos\left(\tfrac{\pi}{2}(a+b)\right) \left[ I(a,b,c)+J(a,b,c) \right]\Big].
\end{multline}

\end{Thm}
\begin{proof}
The difference of the two expressions stated in Theorem \ref{thm31} yield
\begin{align}
\label{T-explicit-2}
2 f_{c}(c) T(a,b,c) =
&\left( \frac{1}{f_{c}(a)f_{c}(b)} - \frac{1}{f_{s}(a)f_{s}(b)} \right)\left[ I(a,b,c) + J(a,b,c) \right] \\
\nn
+ &\left( \frac{1}{f_{c}(a)f_{c}(b)} + \frac{1}{f_{s}(a)f_{s}(b)} \right)\left[ J(c,a,b) + J(c,b,a) \right]
\end{align}
\no
and the result follows directly from here. The values of
$a, b, c \in \mathbb{N}$
are excluded due to the singularity of $\lambda(z)$ for $z \in \mathbb{N}$.
\end{proof}

\section{The limiting case} \label{S:limit}

The goal of this section is to analyze the result of Theorem \ref{thm36}
as the parameters $a, \, b, \, c$ approach positive integer values.  The
notation
$a = n_{1} + \varepsilon_{1}, \;
b = n_{2} + \varepsilon_{2}, \;
c = n_{3} + \varepsilon_{3}$
\no
with $n_{j} \in \mathbb{N}$ and $\varepsilon_{j} \to 0$ is used. \\

We start by writing
\ba\label{rel-T-T-tilde}
T(a,b,c): = \frac{{(2\pi )^{a + b + c} }}
{{16\Gamma (a)\Gamma (b)\Gamma (c)}}\tilde T(a,b,c)
\ea
with, according to \eqref{deffcfs} and \eqref{T-explicit-2},
\begin{multline}\label{T-tilde}
\tilde T(a,b,c): = \frac{1}
{{c_c }}\left[ {\frac{{I(a,b,c) + J(a,b,c) + J(c,a,b) + J(c,b,a)}}
{{c_a c_b }}} \right. \\
+ \left. {\frac{{ - I(a,b,c) - J(a,b,c) + J(c,a,b) + J(c,b,a)}}
{{s_a s_b }}} \right]
\end{multline}

\no
where
\[
c_a  = \cos (\pi a/2)\quad{\text{and}}\quad s_a  = \sin (\pi a/2)
\]
\no
and $c_{b}, s_{b}, c_{c}, s_{c}$ are similarly defined.  \\

Our first task will be to obtain the limit of $c_c \tilde T(a,b,c)$ as
both $a$ and $b$ approach positive integer values.
The functions $c_a$ and $s_a$ have the property that, as the real number
$a$ approaches an integer, one of them tends to zero while the other tends to
$+1$ or $-1$, depending on the parity.
Explicitly, for $n$ an integer and $\varepsilon$ an infinitesimal
quantity,
\ba
\label{limit-ca}
c_{n+\varepsilon} & = &
\begin{cases}
(-1)^{n/2}+o(\varepsilon),\quad n \text{  even} \\
(-1)^{\frac{n+1}{2}}\displaystyle{\frac{\pi}{2}}\varepsilon+o(\varepsilon),\quad n \text{  odd},
\end{cases}
\\
\label{limit-sa}
s_{n+\varepsilon} & = &
\begin{cases}
(-1)^{n/2}\displaystyle{\frac{\pi}{2}}\varepsilon+o(\varepsilon),\quad n \text{  even} \\
(-1)^{\frac{n-1}{2}}+o(\varepsilon),\quad n \text{  odd},
\end{cases}
\ea
\no
Hence, in order to compute the limit we are seeking, we need to expand the numerators
inside square brackets in \eqref{T-tilde} up to order $\varepsilon_1\varepsilon_2$.
This is accomplished by replacing both $\barzeta(a,q)$ and $\barzeta(b,q)$ by their
Taylor series expansions around an integer value of its first argument,
\ba
\nn
\barzeta (n + \varepsilon ,q) &=& \barzeta (n,q) + \varepsilon \barzeta '(n,q) + o(\varepsilon )\\
 &=&  - \frac{1}{n}\left[ {B_n (q) + \varepsilon A_n (q)} \right] + o(\varepsilon ),
\ea
according to \eqref{Hurwitz-to-Bernoulli} and the definition of the Bernoulli function
\ba\label{def-A}
A_k (q): = k\zeta '(1 - k,q),
\ea
studied in \cite{esmo2,esmo4}.
For instance, using the shorthand notation
\ba\label{integration}
\left\langle {f(q)} \right\rangle : = \int_0^1 {f(q)\,dq}
\ea
we have
\begin{align}
J(c,a,b) \Big|_{\small{\begin{matrix}a = n_1  + \varepsilon _1 \\
b = n_{2} + \varepsilon _2 \end{matrix}}}  = &
\frac{1}{{n_1 n_2 }}\left[
{\left\langle {B_{n_1 } (q)B_{n_2 } (1 - q)\barzeta \left( {c,q} \right)} \right\rangle } \right. \\
\nn
&\qquad\quad+ \,\varepsilon _1 \left\langle {A_{n_1 } (q)B_{n_2 } (1 - q)\barzeta \left( {c,q} \right)} \right\rangle\\
\nn
&\qquad\quad+ \,\varepsilon _2 \left\langle {B_{n_1 } (q)A_{n_2 } (1 - q)\barzeta \left( {c,q} \right)} \right\rangle\\
\nn
&\qquad\quad+ \,\varepsilon _1 \varepsilon _2
\left. \left\langle {A_{n_1 } (q)A_{n_2 } (1 - q)\barzeta \left( {c,q} \right)} \right\rangle \right]
+o(\varepsilon _1 \varepsilon _2).
\end{align}
Using the reflection property \eqref{Bernoulli-reflection} of the Bernoulli polynomials
and the invariance of the integration \eqref{integration} under the change of
variable $q\to 1-q$, we find the following result for the numerators inside square
brackets in \eqref{T-tilde} (the upper sign corresponds to the numerator of the
first term and the lower sign corresponds to the numerator of the second term):
\begin{multline}\label{numerators}
\pm I(a,b,c) \pm J(a,b,c) + J(c,a,b) + J(c,b,a) = \\
\frac{1}{{n_1 n_2 }}
\Big\{
\pm \left[ { \pm 1 + ( - 1)^{n_1 } } \right]\left[ { \pm 1 + ( - 1)^{n_2 } } \right]
\left\langle {B_{n_1 } (q)B_{n_2 } (q)\barzeta \left( {c,q} \right)} \right\rangle
\\
+ \varepsilon _1 \left[ { \pm 1 + ( - 1)^{n_2 } } \right]
\left\langle {A_{n_1 } (q)B_{n_2 } (q)\barzeta _ +  \left( {c,q} \right)} \right\rangle
\\
+ \varepsilon _2 \left[ { \pm 1 + ( - 1)^{n_1 } } \right]
\left\langle {B_{n_1 } (q)A_{n_2 } (q)\barzeta _ +  \left( {c,q} \right)} \right\rangle
\\
+ \varepsilon _1 \varepsilon _2
\left[
{ \pm \left\langle {A_{n_1 } (q)A_{n_2 } (q)\barzeta _ +  \left( {c,q} \right)} \right\rangle }
+
{\left\langle {A_{n_1 } (q)A_{n_2 } (1 - q)\barzeta _ +  \left( {c,q} \right)} \right\rangle } \right]
\Big\}
\\
+o(\varepsilon _1 \varepsilon _2).
\end{multline}
where
\ba\label{def-barzeta-plus}
\barzeta_+  (z,q): = \barzeta \left( {z,q} \right) + \barzeta \left( {z,1 - q} \right).
\ea

We now examine the behavior of the Tornheim sum as
$\varepsilon_1, \varepsilon_2 \to 0$. The limiting value is obtained from
\eqref{limit-ca}, \eqref{limit-sa} and \eqref{numerators}. Observe
that $T(a,b,c) = T(b,a,c)$ so only three cases are
presented. \\

\begin{Thm}
Suppose $n_{1}, \, n_{2} \in \mathbb{N}$ and
$c \in \mathbb{R} \backslash \mathbb{N}$. Then the
Tornheim double series $T(n_{1},n_{2},c)$ are given by \\

\begin{multline}
( - 1)^{\frac{{n_1  + n_2 }}{2}} \frac{{(2\pi )^{n_1  + n_2 } }}{{4n_1 !n_2 !}}
\frac{{(2\pi )^c }}{{\Gamma (c)\cos (\pi c/2)}}
\left[ {\int_0^1 {B_{n_1 } (q)B_{n_2 } (q)\barzeta \left( {c,q} \right)} \,dq} \right.\\
- \frac{1}{{\pi ^2 }}\int_0^1 {A_{n_1 } (q)A_{n_2 } (q)\barzeta_+  \left( {c,q} \right)} \,dq
+ \frac{1}{{\pi ^2 }}\left. {\int_0^1 {A_{n_1 } (q)A_{n_2 } (1 - q)\barzeta_+  \left( {c,q} \right)} \,dq} \right]
\end{multline}
\no
for $n_{1}, \, n_{2}$ even;

\begin{multline}
( - 1)^{\frac{{n_1  + n_2  + 1}}{2}} \frac{{(2\pi )^{n_1  + n_2 } }}{{4n_1 !n_2 !}}
\frac{{(2\pi )^c }}{{\Gamma (c)\cos (\pi c/2)}}
\left[ {\frac{1}{\pi }\int_0^1 {B_{n_1 } (q)A_{n_2 } (q)\barzeta_+  \left( {c,q} \right)} \,dq} \right.\\
+ \frac{1}{\pi }\left. {\int_0^1 {A_{n_1 } (q)B_{n_2 } (q)\barzeta_+  \left( {c,q} \right)} \,dq} \right]
\end{multline}
\no
for $n_{1}$ even and $n_{2}$ odd, and

\begin{multline}
( - 1)^{\frac{{n_1  + n_2 }}{2}} \frac{{(2\pi )^{n_1  + n_2 } }}{{4n_1 !n_2 !}}
\frac{{(2\pi )^c }}{{\Gamma (c)\cos (\pi c/2)}}
\left[ {\int_0^1 {B_{n_1 } (q)B_{n_2 } (q)\barzeta \left( {c,q} \right)} \,dq} \right.\\
- \frac{1}{{\pi ^2 }}\int_0^1 {A_{n_1 } (q)A_{n_2 } (q)\barzeta_+  \left( {c,q} \right)} \,dq
- \frac{1}{{\pi ^2 }}\left. {\int_0^1 {A_{n_1 } (q)A_{n_2 } (1 - q)\barzeta_+  \left( {c,q} \right)} \,dq} \right]
\end{multline}
\no
for $n_{1}, \, n_{2}$ odd;
\end{Thm}

\medskip

The final step in the process is to let $c = n_{3} + \varepsilon_{3}$ and let
$\varepsilon_{3} \to 0$. For $n$ even we simply have
\[
\frac{{\barzeta_+  (n,q)}}{{\cos (\pi n/2)}} =  - \frac{2}{n}( - 1)^{n/2} B_n (q),
\]
whereas for $n$ odd,
\[
\lim_{c\to n} \frac{{\barzeta_+  (c,q)}}
{{\cos (\pi c/2)}} =  - \frac{2}
{n}( - 1)^{\frac{{n + 1}}
{2}} \frac{1}
{\pi }\left[ {A_n (q) + A_n (1 - q)} \right].
\]

The value of $T(n_{1},n_{2},n_{3})$ is thus
expressed in terms of integrals of triple products of the Bernoulli
polynomials $B_{k}(q) = -k \zeta(1-k,q)$ and the function
$A_{k}(q)  =  k \zeta'(1-k,q)$.   \\

Define the following families of integrals:

{\allowdisplaybreaks
\begin{align}
\label{def-R1}
R_{1}(n_{1},n_{2},n_{3}) & =  \int_{0}^{1} B_{n_{1}}(q) B_{n_{2}}(q) B_{n_{3}}(q) \, dq   \\
\label{def-R2}
R_{2}(n_{1},n_{2},n_{3}) & =  \frac{1}{\pi}\int_{0}^{1} B_{n_{1}}(q) B_{n_{2}}(q) A_{n_{3}}(q) \, dq \\
\label{def-R3}
R_{3}(n_{1},n_{2},n_{3}) & =  \frac{1}{\pi^2}\int_{0}^{1} A_{n_{1}}(q) A_{n_{2}}(q) B_{n_{3}}(q) \, dq \\
\label{def-R4}
R_{4}(n_{1},n_{2},n_{3}) & =  \frac{1}{\pi^2}\int_{0}^{1} A_{n_{1}}(q) A_{n_{2}}(1-q) B_{n_{3}}(q) \, dq \\
\label{def-R5}
R_{5}(n_{1},n_{2},n_{3}) & =  \frac{1}{\pi^3}\int_{0}^{1} A_{n_{1}}(q) A_{n_{2}}(q) A_{n_{3}}(q) \, dq \\
\label{def-R6}
R_{6}(n_{1},n_{2},n_{3}) & =  \frac{1}{\pi^3}\int_{0}^{1} A_{n_{1}}(q) A_{n_{2}}(q) A_{n_{3}}(1-q) \, dq.
\end{align}
}
These integrals are all symmetric under interchange of their first two arguments
($n_1$ and $n_2$), except for $R_4$ which is antisymmetric if $n_3$ is odd.

\no
Define
\ba
p(n) & = &
\begin{cases}
(-1)^{n/2},\quad n \text{  even} \\
(-1)^{\frac{n+1}{2}},\quad n \text{  odd}.
\end{cases}
\ea
\no

\begin{Thm}
\label{thm-zagr}
Let $\alpha = n_{1}+n_{2}+n_{3}$. Then
the Tornheim double series $T(n_{1},n_{2},n_{3})$ is given by
\ba
T(n_1 ,n_2 ,n_3 ) = p(\alpha)\frac{{(2\pi )^{\alpha} }}
{{2n_1!\,n_2!\,n_3!}}T_R (n_1 ,n_2 ,n_3 )
\ea
where $T_R(n_1 ,n_2 ,n_3 )$ can be expressed in terms
of the functions $R_{j}: 1 \leq j \leq 6$ as follows:

{\allowdisplaybreaks
\begin{align}
\intertext{{\bf Case 1}. $n_{1}, \, n_{2}$ and $n_{3}$ are even:}
T_R (n_1 ,n_2 ,n_3 ) &=  - \frac{1}
{2}R_1 (n_1 ,n_2 ,n_3 ) + R_3 (n_1 ,n_2 ,n_3 ) - R_4 (n_1 ,n_2 ,n_3 )\\
\intertext{{\bf Case 2}. $n_{1}$ and $n_{2}$ are even; $n_{3}$ is odd:}
T_R (n_1 ,n_2 ,n_3 ) &= - R_2 (n_1 ,n_2 ,n_3 ) + R_5 (n_1 ,n_2 ,n_3 )\\
&\quad+ R_6 (n_1 ,n_2 ,n_3 ) - R_6 (n_3 ,n_1 ,n_2 ) - R_6 (n_3 ,n_2 ,n_1 )\nn\\
\intertext{{\bf Case 3}. $n_{1}$ is even, $n_{2}$ is odd, and $n_{3}$ is even:}
T_R (n_1 ,n_2 ,n_3 ) &=  - R_2 (n_3 ,n_1 ,n_2 ) - R_2 (n_3 ,n_2 ,n_1 )\\
\intertext{{\bf Case 4}. $n_{1}$ is even; $n_{2}$ and $n_{3}$ are odd:}
T_R (n_1 ,n_2 ,n_3 ) &= R_3 (n_3 ,n_1 ,n_2 ) + R_3 (n_3 ,n_2 ,n_1 )\\
&\quad+ R_4 (n_1 ,n_3 ,n_2 ) + R_4 (n_2 ,n_3 ,n_1 )\nn\\
\intertext{{\bf Case 5}. $n_{1}$ and $n_{2}$ are odd; $n_{3}$ is even:}
T_R (n_1 ,n_2 ,n_3 ) &=  - \frac{1}
{2}R_1 (n_1 ,n_2 ,n_3 ) + R_3 (n_1 ,n_2 ,n_3 ) + R_4 (n_1 ,n_2 ,n_3 )\\
\intertext{{\bf Case 6}. $n_{1}, \, n_{2}$ and $n_{3}$ are odd:}
T_R (n_1 ,n_2 ,n_3 ) &=  - R_2 (n_1 ,n_2 ,n_3 ) + R_5 (n_1 ,n_2 ,n_3 )\\
&\quad+ R_6 (n_1 ,n_2 ,n_3 ) + R_6 (n_3 ,n_1 ,n_2 ) + R_6 (n_3 ,n_2 ,n_1 ). \nn
\end{align}
} 
\end{Thm}

\medskip

The closed form evaluation of the
Tornheim sums $T(n_{1},n_{2},n_{3})$ has thus been
reduced to that of the integrals $R_{j}$. A partial evaluation
of these integrals is presented in the next section, in terms of
new family  of integrals, closely related to $R_{j}$.

\section{A new family of integrals} \label{S:Q-integrals}

The evaluation of the integrals $R_j$ is most conveniently organized
in terms of a new family of integrals $Q_{j}$, defined in terms of
the balanced generalized polygamma function, introduced in \cite{esmo4},
by
{\allowdisplaybreaks
\begin{align}
\label{def-Q1}
Q_{1}(n_{1},n_{2},n_{3}) & = \int_{0}^{1} B_{n_{1}}(q) B_{n_{2}}(q)
B_{n_{3}}(q) \, dq  \\
\label{def-Q2}
Q_{2}(n_{1},n_{2},n_{3}) & = \int_{0}^{1} B_{n_{1}}(q) B_{n_{2}}(q)
\psi^{(-n_{3})}(q)\, dq \\
\label{def-Q3}
Q_{3}(n_{1},n_{2},n_{3}) & = \int_{0}^{1} \psi^{(-n_{1})}(q)
\psi^{(-n_{2})}(q) B_{n_{3}}(q) \, dq \\
\label{def-Q4}
Q_{4}(n_{1},n_{2},n_{3}) & = \int_{0}^{1} \psi^{(-n_{1})}(q)
\psi^{(-n_{2})}(1-q) B_{n_{3}}(q) \, dq  \\
\label{def-Q5}
Q_{5}(n_{1},n_{2},n_{3}) & = \int_{0}^{1} \psi^{(-n_{1})}(q)
\psi^{(-n_{2})}(q) \psi^{(-n_{3})}(q) \, dq  \\
\label{def-Q6}
Q_{6}(n_{1},n_{2},n_{3}) & = \int_{0}^{1} \psi^{(-n_{1})}(q)
\psi^{(-n_{2})}(q) \psi^{(-n_{3})}(1-q) \, dq.
\end{align}
} 
\no
These integrals satisfy the same symmetry properties as their $R$-analogs.
In addition, the $Q$-integrals satisfy homogeneous recursion relations, which allow
their evaluation in terms of a few basic ones.
The relation among the families
$R_{j}$ and $Q_{j}$ using the identity
\ba
A_{m}(q) & = & m! \psi^{(-m)}(q) + h_{m}B_{m}(q)  \label{relation1}
\ea
\no
is given in Appendix 3.

In this section we present recurrences for the integrals $Q_{j}$. The initial
conditions require a variety of definite integrals listed below:
\begin{align}
\intertext{$\bullet$  The integrals}
\label{def-N}
N_{m,n}  & = \int_{0}^{1}  \psi^{(-m)}(q)B_{n}(q) \, dq,\\
\label{def-M}
M_{m,n}  & =  \int_{0}^{1} \psi^{(-m)}(q) \psi^{(-n)}(q) \, dq,  \\
\intertext{and}
\label{def-M*}
M^{*}_{m,n}  & =  \int_{0}^{1} \psi^{(-m)}(q) \psi^{(-n)}(1-q) \, dq,\\
\intertext{which will be evaluated in Section \ref{sec-eval-zetapsi}.}
\intertext{$\bullet$ The families of integrals}
\label{def-K}
K_{m,n} & = \ione  \psi^{(-m)}(q)B_{n}(q) \,\ln \Gamma(q)  \, dq,  \\
\label{def-K*}
K^{*}_{m,n} & = \ione  \psi^{(-m)}(1-q) B_{n}(q)\,\ln \Gamma(q)  \, dq,\\
\label{def-Z}
Z_{m,n} & = \ione \psi^{(-m)}(q)  \psi^{(-n)}(q) \ln \Gamma(q) \, dq,  \\
\intertext{and}
\label{def-Z*}
Z_{m,n}^{*} & = \ione \psi^{(-m)}(q)  \psi^{(-n)}(1-q) \ln \Gamma(q) \, dq.
\end{align}
The closed form evaluation of these functions is left as an open question.

\bigskip

\no
{\bf The integral} $Q_{1}$. The
explicit value of
$Q_{1}(n_{1},n_{2},n_{3})$ was given by Carlitz \cite{carlitz}:
\begin{align}
\label{Q1-explicit}
Q_{1}(n_{1},n_{2},n_{3})
& =
(-1)^{n_{3}+1} n_{3}!
\sum_{k=0}^{\lfloor{ (n_{1}+n_{2}-1)/2 \rfloor} }
\left[  n_{1} \binom{n_{2}}{2k} +
n_{2} \binom{n_{1}}{2k} \right] \\
& \hspace{0.5in} \times \frac{(n_{1}+n_{2}-2k-1)!}{(n_{1}+n_{2}+n_{3}-2k)!}
B_{2k}B_{n_{1}+n_{2}+n_{3}-2k}. \nn
\end{align}

\medskip

\no
{\bf The integral}  $Q_{2}$.  This integral is obtained directly from the formula
given in \cite{esmo4},
\begin{multline}
\label{psi-int-3}
\int_0^1 \zeta (z+1,q)\psi (z',q)dq = 2(2\pi )^{z + z'}\Gamma (-z)
\Bigg\{
\frac{\pi }{2}\zeta (- z - z')\sin \frac{\pi }{2}(z - z')\\
+
\Big[ (\gamma  + \ln 2\pi )\zeta (-z - z')
- \zeta '(- z - z') \Big]\cos \frac{\pi }{2}(z - z')
\Bigg\},
\end{multline}
valid for $\realpart z, \realpart z'<0$ and $\realpart (z+z')<-1$.
The evaluation at $z=-m$ and $z'=-n$, with $m,n\in\allN$ gives
\begin{multline}
\label{psi-int-4}
\int_0^1 {B_m (q)\psi^{( - n)}(q)} \,dq = \frac{{2m!}}{{(2\pi )^{m + n} }}
\left\{ {\frac{\pi }{2}} \zeta \left( {m + n} \right)\sin \frac{\pi }{2}(m - n)\right.\\
- \left[ {(\gamma  + \ln 2\pi )\zeta \left( {m + n} \right) - \zeta '\left( {m + n} \right)} \right]
\left. {\cos \frac{\pi }{2}(m - n)} \right\}.
\end{multline}

The evaluation of $Q_2(n_{1},n_{2},n_{3})$ follows from \eqref{psi-int-4} and the
representation \eqref{prod-2-Bernoullis} for the product of two Bernoulli
polynomials:

\medskip

\begin{Thm}
Let $n_{1}, n_{2}, n_{3} \in \mathbb{N}$ and let
$\alpha = n_{1}+n_{2}+n_{3}$. Then
\begin{multline}
\label{Q2-explicit}
Q_2 (n_1 ,n_2 ,n_3 ) = 2(-1)^{n_3}\sum\limits_{k = 0}^{k(n_1 ,n_2 )}
\left[  n_{1} \binom{n_{2}}{2k} + n_{2} \binom{n_{1}}{2k} \right]
\frac{{(n_1  + n_2  - 2k - 1)!}}{{(2\pi )^{\alpha  - 2k} }}\\
\times
(-1)^k B_{2k}
\left\{ {\frac{\pi }{2}} \sin \frac{{\pi \alpha }}{2}\zeta \left( {\alpha  - 2k} \right)\right.
- \left. {\cos \frac{{\pi \alpha }}{2}
\left[ {(\gamma  + \ln 2\pi )\zeta \left( {\alpha  - 2k} \right) - \zeta '\left( {\alpha  - 2k} \right)} \right]} \right\}
\end{multline}
\no
where $k(n_{1},n_{2}) = \text{Max}\{ \lfloor{ n_{1}/2
\rfloor}, \lfloor{n_{2}/2 \rfloor} \}$.
The constant term in \eqref{prod-2-Bernoullis} gives no contribution on account
that the function $\psi^{(-n)}(q)$ is balanced for $n\in\allN$.
The integrals $Q_2$ are only needed in the case that $\alpha=n_1+n_2+n_3$ is odd,
equal to $2N+1$, say.
In this case, $\sin(\pi\alpha/2)=(-1)^N$ and $\cos(\pi\alpha/2)=0$, so that
\begin{multline}
\label{Q2-explicit-alpha-odd}
Q_2 (n_1 ,n_2 ,n_3 ) = \pi(-1)^{N+n_3}\\
\times\sum\limits_{k = 0}^{k(n_1 ,n_2 )}
\left[  n_{1} \binom{n_{2}}{2k} + n_{2} \binom{n_{1}}{2k} \right]
\frac{{(n_1  + n_2  - 2k - 1)!}}{{(2\pi )^{\alpha  - 2k} }}
(-1)^k B_{2k}\zeta \left( {\alpha  - 2k} \right).
\end{multline}
\no
\end{Thm}
\begin{proof}
The details are elementary.
\end{proof}

\medskip

\no
{\bf The integral}  $Q_{3}$. We now produce a recurrence for
\ba
Q_{3}(n_{1},n_{2},n_{3}) & = & \int_{0}^{1} \psi^{(-n_{1})}(q)
\psi^{(-n_{2})}(q) B_{n_{3}}(q) \, dq. \nn
\ea
\no
The basic tools are the relations
\begin{align}
\frac{d}{dq} \psi^{(-m)}(q) & = \psi^{(-m+1)}(q),\\
\intertext{and}
\frac{d}{dq} B_m(q)&= m B_{m-1}(q),
\end{align}
\no
valid for $m \in \nnN$,
and the fact that both the negapolygamma functions and the Bernoulli
polynomials are balanced for the range of indices we wish to consider, i.e.,
$\psi^{(-m)}(1) = \psi^{(-m)}(0)$, for $m\ge 2$, and
$B_m(1)=B_m(0)$, for all $m$.\\

\begin{Thm}
Let $n_{1}, \, n_{2}, \, n_{3} \in \mathbb{N}$ with $n_{1}, n_{2} > 1$.
Then
\begin{multline}
\label{recq3}
(n_{3}+1)Q_{3}(n_{1},n_{2},n_{3}) = -Q_{3}(n_{1}-1,n_{2},n_{3}+1)
-Q_{3}(n_{1},n_{2}-1,n_{3}+1).
\end{multline}
\end{Thm}
\begin{proof}
Start with
\ba
(n_{3}+1)Q_{3}(n_{1},n_{2},n_{3}) & =  & \ione  \psi^{(-n_{1})}(q)
\, \psi^{(-n_{2})}(q) \frac{d}{dq}B_{n_{3}+1}(q) \, dq, \nn
\ea
\no
integrate by parts and observe that there is no contribution from the
boundary.
\end{proof}

\medskip

The recurrence shows that the value of $Q_{3}(n_{1},\,n_{2},n_{3})$
can be obtained from the values of
\ba
\label{initcond-q3-2}
Q_{3}(1,m,n) & = &
K_{m,n} + \zeta'(0) N_{m,n},
\ea
\no
in view of $\psi^{(-1)}(q) = \ln \Gamma(q) + \zeta'(0)$, the
symmetry of the integral $Q_3$ under interchange of its first two arguments,
and the definitions \eqref{def-K} and \eqref{def-N} of the integrals
$K_{m,n}, N_{m,n}$.

\medskip

\no
{\bf The integral}  $Q_{4}$. Similarly, for $n_{1}, \, n_{2} > 1$,
we have that
\ba
Q_{4}(n_{1},n_{2},n_{3}) & = & \int_{0}^{1} \psi^{(-n_{1})}(q)
\psi^{(-n_{2})}(1-q) B_{n_{3}}(q) \, dq  \nn
\ea
\no
satisfies the recurrence
\begin{multline}
\label{recq4}
(n_{3}+1)Q_{4}(n_{1},n_{2},n_{3}) = -Q_{4}(n_{1}-1,n_{2},n_{3}+1)
+Q_{4}(n_{1},n_{2}-1,n_{3}+1),
\end{multline}
\no
so that it can be obtained from
\ba
\label{initcond-q4-2}
Q_{4}(1,m,n)
& = & K^{*}_{m,n} + \zeta'(0) (-1)^n N_{m,n},
\ea
where $K^{*}_{m,n}$ is defined in \eqref{def-K*}.

\medskip

\no
{\bf The integrals} $Q_{5}$ and $Q_{6}$. Similar arguments show that for
$n_{1}, \, n_{2} > 1$ we have
\ba
Q_{5}(n_{1},n_{2},n_{3}) & = & \int_{0}^{1} \psi^{(-n_{1})}(q)
\psi^{(-n_{2})}(q) \psi^{(-n_{3})}(q) \, dq  \nn \\
Q_{6}(n_{1},n_{2},n_{3}) & = & \int_{0}^{1} \psi^{(-n_{1})}(q)
\psi^{(-n_{2})}(q) \psi^{(-n_{3})}(1-q) \, dq  \nn
\ea
\no
satisfy the recurrences
\begin{align}
\label{recq5}
Q_{5}(n_{1},n_{2},n_{3}) & = - Q_{5}(n_{1}-1,n_{2},n_{3}+1) -
Q_{5}(n_{1},n_{2}-1,n_{3}+1)\\
\intertext{and}
\label{recq6}
Q_{6}(n_{1},n_{2},n_{3}) & = Q_{6}(n_{1}-1,n_{2},n_{3}+1) +
Q_{6}(n_{1},n_{2}-1,n_{3}+1).
\end{align}
\no

\medskip

The initial conditions are
\ba
Q_{5}(1,m,n) & = &
 Z_{m,n} + \zeta'(0) M_{m,n}, \\
Q_{6}(1,m,n) & = &
 Z^{*}_{m,n} + \zeta'(0) M^{*}_{m,n},
\ea
where $Z_{m,n}$, $Z^{*}_{m,n}$, $M_{m,n}$ and $M^{*}_{m,n}$ are defined
in \eqref{def-Z}, \eqref{def-Z*}, \eqref{def-M} and \eqref{def-M*},
respectively.

\section{Evaluation of integrals}
\label{sec-eval-zetapsi}

\subsection{The product of a Bernoulli polynomial and a balanced ne\-ga\-poly\-gam\-ma}
\medskip

The explicit evaluation of the integral $N_{m,n}$ defined by \eqref{def-N} was
obtained in \eqref{psi-int-4} as a byproduct of the evaluation of $Q_2$:
\begin{multline}
\label{eval-N}
N_{m,n}= \frac{{2n!}}{{(2\pi )^{m + n} }}
\left\{ {\frac{\pi }{2}} \zeta \left( {m + n} \right)\sin \frac{\pi }{2}(n - m)\right.\\
- \left[ {(\gamma  + \ln 2\pi )\zeta \left( {m + n} \right) - \zeta '\left( {m + n} \right)} \right]
\left. {\cos \frac{\pi }{2}(n - m)} \right\}.
\end{multline}

\subsection{The product of two balanced negapolygammas}
\medskip

The integral of the product of two balanced ne\-ga\-pol\-y\-gam\-ma
functions has been given in \cite{esmo2}:

\no
Let $k, \, k' \in \mathbb{N}, \,  k_{+} = k+k'$  and $k_{-} = k -k'$. Then
\begin{align}
\label{eval-M}
M_{k,k'} & =\ione\psi ^{( - k)} (q)\psi ^{( - k')} (q)\,dq \\
&= \frac{2}{(2 \pi)^{k_{+}}}
\cos \left( \frac{\pi }{2}k_{-} \right)
\Big[ A_{+} \zeta(k_{+}) -2A \zeta'(k_{+}) + \zeta''(k_{+}) \Big] \nn
\intertext{where $A = \gamma + \ln 2 \pi$ and $A_\pm = A^{2}  \pm \pi^{2}/4$.}
\intertext{Similarly, one finds}
\label{eval-M*}
\nn
M^{*}_{k,k'} & = \ione\psi ^{(-k)} (q)\psi ^{(-k')} (1-q)\,dq \stvictor\\
& = \frac{2}{(2 \pi)^{k_{+}}}
\Big\{
\cos \left( \frac{\pi }{2}k_{+} \right)
\Big[ A_{-} \zeta(k_{+}) -2A \zeta'(k_{+}) + \zeta''(k_{+}) \Big]\\
\nn
&\qquad\qquad\qquad + \pi \, \sin \left(\frac{\pi }{2}k_{+} \right)
\left[ A \zeta(k_{+}) - \zeta'(k_{+}) \right]
\Big\}.
\end{align}

\medskip

We have been unable to obtain closed form results for the integrals $K_{m,n}$
$K_{m,n}^*$, $Z_{m,n}$ and $Z_{m,n}^*$, which involve the kernel $\ln\Gamma(q)$
and one or two negapolygamma functions.

\section{Some examples} \label{sec-examples}

In this section we describe some explicit evaluations of Tornheim sums.
It is convenient to introduce the function
\ba
U_{m,n} & = & \int_{0}^{1}  \psi^{(-m)}(q) B_{n}(q)\, \ln \sin( \pi q) \;
dq.
\ea
\no
The identity
\ba
\ln \Gamma (q) + \ln \Gamma(1-q) & = & \ln \pi - \ln \sin (\pi q)
\ea
\no
yields the relation
\ba
K_{m,n} + (-1)^{n} K^{*}_{m,n} & = & \ln \pi N_{m,n} - U_{m,n}.
\ea

\medskip

\begin{Example}
We consider the value of $T(1,1,2)$.
This corresponds to Case 5 of Theorem \ref{thm-zagr}. It yields
\begin{align}
T(1,1,2) & = 4\pi^{4} \left[ -\frac{1}{2}R_{1}(1,1,2)
+R_{3}(1,1,2) + R_{4}(1,1,2) \right],\\
\intertext{and in terms of the $Q_{j}$-family,}
T(1,1,2) & = 4\pi^{4}\left[ -\frac{1}{2}Q_{1}(1,1,2)
+\frac{1}{\pi^2}Q_{3}(1,1,2) + \frac{1}{\pi^2}Q_{4}(1,1,2) \right].
\end{align}
\no
The $Q_{j}$-integrals are given by
\ba
Q_{1}(1,1,2)  &= & \frac{1}{180} \nn \\
Q_{3}(1,1,2) & = & K_{1,2} + \zeta'(0) N_{1,2} \nn \\
Q_{4}(1,1,2) & = & K^{*}_{1,2} + \zeta'(0) N_{1,2} \nn
\ea
\no
and using the values $\zeta'(0) = -\frac{1}{2} \ln (2 \pi)$ and
$N_{1,2} = \zeta(3)/4 \pi^{2}$ we obtain
\ba
T(1,1,2) & = & 4 \pi^{2} \left( K_{1,2} + K^{*}_{1,2} \right) -
\zeta(3) \ln ( 2 \pi) - \frac{1}{90} \pi^{4}.  \nn
\ea
\no
In terms of the $U$-function this can be written as
\ba
T(1,1,2) & = &  - 4 \pi^{2} U_{1,2} -\zeta(3) \ln 2 - \frac{1}{90} \pi^{4}.
\ea

The identities \eqref{elem00} and \eqref{symmetryrule} yield
$T(1,1,2) = 2 T(0,1,3)$, and the method of Huard et al. \cite{huard}
yields the values of
$T(n_{1},n_{2},n_{3})$ for $N = n_{1}+n_{2}+n_{3}$ odd and also for $N=4$
and $6$. For instance $T(0,1,3) = \tfrac{1}{4} \zeta(4)$. This yields an
evaluation of an integral of type
$U_{m,n}$: the value $T(1,1,2) = \zeta(4)/2$, the identity
$\psi^{(-1)}(q) = \ln \Gamma (q) + \zeta'(0)$ and
\ba\nn
\ione B_{2}(q) \, \ln ( \sin \pi q ) \, dq & = & - \frac{\zeta(3)}{2 \pi^{2}}
\ea
\no
given in Example $5.2$ of \cite{esmo1} produce
\ba
\ione B_{2}(q) \ln \Gamma(q) \, \ln ( \sin \pi q ) \; dq & = &
- \left( \frac{\pi^{2}}{240} + \frac{\ln( 4 \pi) \, \zeta(3) }{4 \pi^{2} }
\right).
\ea
\no
It follows that
\ba
U_{1,2} & = & -\frac{\pi^{2}}{240} - \frac{\ln 2 \, \zeta(3) }{4 \pi^{2}}.
\ea
\end{Example}

\medskip

\begin{Example}
The explicit expression for $R_{2}(n_{1}, n_{2},n_{3})$ that can be obtained from
\eqref{rel-RQ2} through \eqref{Q1-explicit} and \eqref{Q2-explicit} permits the evaluation
of $T(n_{1},n_{2},n_{3})$ in the case
$n_{1}, n_{3}$ even and $n_{2}$ odd. For example
\begin{align}
T(2,1,2) & =\frac{\pi^{2}}{6}  \, \zeta(3) - \frac{3}{2} \zeta(5), \nn \\
T(2,3,2) & = - \frac{\pi^{2}}{6}  \zeta(5) + 2 \zeta(7) , \stvictor \nn\\
\intertext{and}
T(4,3,2) & =
\frac{\pi^{4}}{90}  \zeta(5) + \frac{\pi^{2}}{6}  \zeta(7) -
\frac{5}{2} \zeta(9). \stvictor \nn
\end{align}
\end{Example}

\medskip

\begin{Example}
Define $w = a+b+c$ to be the {weight} of the sum $T(a,b,c)$. The results of the
procedure described above for sums of small weight are given below.

{\allowdisplaybreaks
\begin{align*}
\intertext{weight 3:}
T(1,1,1) & =
4Z_{1,1}  + 12Z_{1,1}^*  - \zeta \left( 3 \right)
+ \ln{2\pi}\left( {\frac{{A^2 }}{3} - \frac{{\pi ^2 }}{6}
- \frac{{4A\zeta '\left( 2 \right)}}{{\pi ^2 }}
+ \frac{{2\zeta ''\left( 2 \right)}}{{\pi ^2 }}} \right),\\
\intertext{weight 4:}
T(1,1,2) & =
- 4 \pi^{2} U_{1,2} -\frac{\pi^{4}}{90} - \zeta(3)\ln 2 , \\
T(1,2,1) & =  -4 \pi^{2} \left( U_{1,2} +2 U_{2,1}\right),\\
\intertext{weight 5:}
T(1,1,3) &=
- 8\pi ^2 \left( {K_{1,3}^*  + 2Z_{1,3}  + 2Z_{1,3}^*  + 4Z_{3,1}^* } \right)
- \zeta \left( 5 \right)\\
&\quad + \ln{2\pi}
\left(
{\frac{{\pi ^2 A^2 }}{{45}} - \frac{{\pi ^2 A}}{{30}} - \frac{{\pi ^4 }}{{90}}
- \frac{{4\zeta '\left( 4 \right)A}}{{\pi ^2 }}
+ \frac{{3\zeta '\left( 4 \right)}}{{\pi ^2 }}
+ \frac{{2\zeta ''\left( 4 \right)}}{{\pi ^2 }}} \right),\\
T(1,2,2) & =
\frac{\pi^{2}}{6} \zeta(3) - \frac{3}{2} \zeta(5), \\
T(1,3,1) &=
- 8\pi ^2 \left( {K_{1,3}^*  + 2Z_{3,1}  + 2Z_{1,3}^*  + 4Z_{3,1}^* } \right)
+ \frac{\pi^2}{6}\zeta \left( 3 \right)- 2\zeta \left( 5 \right)\\
&\quad + \ln{2\pi}
\left(
{\frac{{\pi ^2 A^2 }}{{45}} - \frac{{\pi ^2 A}}{{30}} - \frac{{\pi ^4 }}{{90}}
- \frac{{4\zeta '\left( 4 \right)A}}{{\pi ^2 }}
+ \frac{{3\zeta '\left( 4 \right)}}{{\pi ^2 }}
+ \frac{{2\zeta ''\left( 4 \right)}}{{\pi ^2 }}} \right),\\
T(2,2,1) &=
32\pi ^2 Z_{2,2}  + \frac{{\pi ^2 }}{3}\zeta \left( 3 \right) - 3\zeta \left( 5 \right)\\
&\quad + \ln{2\pi}
\left(
{ - \frac{{\pi ^2 A^2 }}{{45}} - \frac{{\pi ^4 }}{{180}}
+ \frac{{4\zeta '\left( 4 \right)A}}{{\pi ^2 }}
- \frac{{2\zeta ''\left( 4 \right)}}{{\pi ^2 }}}
\right).
\end{align*}
}

\end{Example}

\medskip

We have produced some partial results in the evaluation of the
integrals $K_{m,n}$, $K^{*}_{m,n}$ and
$Z_{m,n}, \, Z^{*}_{m,n}$. These suggest that the value
of the Tornheim sums can be expressed in terms of  a small number of
definite integrals. For instance, for $m \geq 3$ odd, we have the relation
\ba
K_{m,n} & = & -mK_{m-1,n+1} + \ln \sqrt{2 \pi} ( N_{m,n} + m N_{m-1,n+1} )
\nn \\
 & & - \int_{0}^{1} B_{m}(q) \psi(q) \psi^{(-n-1)}(q) \, dq, \nn
\ea
\no
which reduces the value of $K_{m,n}$ to that of $K_{1,m+n-1}$ plus the moments
of the product $\psi(q) \psi^{(-j)}(q)$. Details
will be presented elsewhere.

\bigskip

\section{Appendix 1: The Bernoulli polynomials} \label{S:berpoly}

The {\em Bernoulli polynomials}
 $B_{n}(q)$ defined by the generating function
\ba
\frac{xe^{xq}}{e^{x}-1} & = & \sum_{n=0}^{\infty} B_{n}(q) \frac{x^{n}}{n!}.
\label{bergen1}
\ea
\no
The {\em Bernoulli numbers} $B_{n} = B_{n}(0)$ satisfy
\ba
B_{n}(q) & = & \sum_{k=0}^{n} \binom{n}{k} B_{k}q^{n-k}.
\label{bergen3}
\ea
\no
For $n \geq 1$ we have the differential recursion
$B_{n}'(q) =  n B_{n-1}(q)$ and the symmetry rule
$B_{n}(1-q)  =  (-1)^{n}B_{n}(q)$. In
particular, $B_{n}(1) = B_{n}(0)$ for $n > 1$.

The Bernoulli polynomials $\{ B_{0}(q), \, B_{1}(q), \, \cdots, B_{n}(q) \}$
form a
basis for the space of polynomials of degree at most $n$. Thus the product
$B_{n_{1}}(q) B_{n_{2}}(q)$ is a linear combination of $B_{j}(q)$ for
$j=0, \cdots, n_{1}+n_{2}$. It is a remarkable fact that this combination has
the explicit form
\begin{multline}
\label{prod-2-Bernoullis}
B_{n_{1}}(q) B_{n_{2}}(q) = \sum_{k=0}^{k(n_{1},n_{2})}
\left[  n_{1} \binom{n_{2}}{2k} +
n_{2} \binom{n_{1}}{2k} \right]
\frac{B_{2k} }{n_{1}+n_{2}-2k} B_{n_{1}+n_{2}-2k}(q)    \\
+ (-1)^{n_{1}+1} \frac{n_{1}! \, n_{2}!}{(n_{1} + n_{2})!} B_{n_{1}+n_{2}},
\end{multline}
\no
where $k(n_{1},n_{2}) = \text{Max}\{ \lfloor{ n_{1}/2
\rfloor}, \lfloor{n_{2}/2 \rfloor} \}$.
In terms of rescaled Bernoulli polynomials and numbers, defined by
\ba
\tilde B_{n}(q)=\frac{B_{n}(q)}{n!},\quad \tilde B_{n}=\tilde B_{n}(0)=\frac{B_{n}}{n!},
\ea
relation \eqref{prod-2-Bernoullis} has the simpler form
\begin{multline}
\label{prod-2-Bernoullis-rescaled}
\tilde B_{n_{1}}(q) \tilde B_{n_{2}}(q) = \sum_{k=0}^{k(n_{1},n_{2})}
\left[  {\textstyle\binom{n_{1}+n_{2}-2k-1}{n_{1}-1}} +
{\textstyle\binom{n_{1}+n_{2}-2k-1}{n_{2}-1}} \right]
\tilde B_{2k}  \tilde B_{n_{1}+n_{2}-2k}(q)    \\
+ (-1)^{n_{1}+1} \tilde B_{n_{1}+n_{2}}.
\end{multline}
\no

In theory, (\ref{prod-2-Bernoullis}) yields expressions for a product of any number of
Bernoulli polynomials. For example,
\ba
\label{square-Bernoulli}
B_{n}^{2}(q) & = & \sum_{k=0}^{\lfloor{n/2 \rfloor}}
\frac{n \binom{n}{2k}}{n-k} B_{2k} B_{2n-2k}(q) +
(-1)^{n+1}
\frac{B_{2n}}{\binom{2n}{n}},
\ea
\no
or
\ba
\label{square-Bernoulli-rescaled}
\tilde B_{n}^{2}(q) & = & 2 \sum_{k=0}^{\lfloor{n/2 \rfloor}}
{\textstyle\binom{2n-2k-1}{n-1}}\tilde B_{2k} \tilde B_{2n-2k}(q) +
(-1)^{n+1}\tilde B_{2n},
\ea
\no
and
\begin{multline}
\label{cube-Bernoulli}
B_{n}^{3}(q) = \sum_{k=0}^{\lfloor{n/2 \rfloor}}
\frac{n \binom{n}{2k}}{n-k}  B_{2k} \;
\sum_{j=0}^{n-k} \left[ n \binom{2n-2k}{2j} + 2(n-k) \binom{n}{2j}\right]
\frac{B_{2j} \, B_{3n-2k-2j}(q) }{3n-2k-2j} \\
+ \;
(-1)^{n+1} \left[
\frac{B_{2n}}{\binom{2n}{n}}B_{n}(q)
+ 2n!^3 \sum_{k=0}^{\lfloor{n/2 \rfloor}}
\binom{2n-2k-1}{n-1}\frac{B_{2k}B_{3n-2k}}{(2k)!(3n-2k)!}
\right].
\end{multline}
or
\begin{align}
\label{cube-Bernoulli-rescaled}
\tilde B_{n}^{3}(q) = 2\sum_{k=0}^{\lfloor{n/2 \rfloor}}
&{\textstyle\binom{2n-2k-1}{n-1}} \tilde B_{2k} \;\\
&\times\sum_{j=0}^{n-k} \left[ {\textstyle\binom{3n-2k-2j-1}{n-1}}
+ {\textstyle\binom{3n-2k-2j-1}{2n-2k-1}}\right]
\tilde B_{2j} \, \tilde B_{3n-2k-2j}(q) \nn\\
&+ \;(-1)^{n+1} \left[
\tilde B_{2n}\tilde B_{n}(q)
+ 2 \sum_{k=0}^{\lfloor{n/2 \rfloor}}
{\textstyle\binom{2n-2k-1}{n-1}}\tilde B_{2k}\tilde B_{3n-2k}
\right].
\nn
\end{align}
\no
Integrating the relation (\ref{bergen1}) yields
\ba
\ione B_{n}(q) \, dq & = & 0, \text{ for } n \geq 1.
\label{evaloneber}
\ea
\no
Apostol \cite{apostol76} gives a direct proof of
\ba
\int_{0}^{1} B_{n_{1}}(q) B_{n_{2}}(q) \, dq & = &
(-1)^{n_{1}+1} \frac{n_{1}! \, n_{2}! }{(n_{1}+n_{2})!} B_{n_{1}+n_{2}},
\ea
for $n_{1}, \, n_{2} \in \allN$.

The Bernoulli polynomials appear also as special values of the
Hurwitz zeta function
\ba
\zeta(1-k,q) &  =& - \frac{1}{k} B_{k}(q).
\label{zetaber}
\ea
\medskip

\section{Appendix 2: The generalized polygamma function $\psi (z,q)$}
\label{S:integrals}

The {\em poly\-gam\-ma function} is defined by
\begin{align}
\label{polygamma-def}
\psi^{(m)}(q) &= \frac{d^{m}}{dq^{m}}
\psi(q), \quad m \in \allN,\\
\intertext{where}
\label{polygamma-hurwitz}
\psi(q) &= \psi^{(0)}(q) = \frac{d}{dq} \ln \Gamma(q)
\end{align}
\no
is the {\em digamma} function.

The function $\psi^{(m)}$ is analytic
in the complex $q$-plane, except for poles (of order
$m+1$) at all non-positive integers. Extensions
of this function for  $m$ a negative integer
have been defined by several authors \cite{adamchik, esmo2, gosper}. These
are the {\em negapolygamma} functions. For example,
Gosper \cite{gosper} defined
\begin{equation}\label{gosper-negapoly}
\begin{split}
\psi_{-1}(q) & = \ln \Gamma(q), \\
\psi_{-k}(q) & = \int_0^q \psi_{-k+1}(t) dt,\quad k\ge 2,\\
\end{split}
\end{equation}
which were later reconsidered by Adamchik \cite{adamchik} in the form
\ba\label{adam-negapoly}
\psi_{-k}(q) = \frac{1}{(k-2)!}\int_0^q (q-t)^{k-2}\ln\Gamma(t) dt,\quad k\ge
2.
\ea
These extensions can be expressed in terms of
the derivative (with respect to its first argument) of the Hurwitz
zeta function at the negative integers
\cite{adamchik,gosper}.
\no
The definition
(\ref{gosper-negapoly}) can be modified by introducing
arbitrary constants of integration at every step. This yields
different extensions differing by polynomials:
\[
\psi _a^{(-m)} (q) - \psi _b^{(-m)} (q) = p_{m - 1} (q),
\]
satisfying
\[
p_n (q) = \frac{d}{{dq}}p_{n + 1} (q).
\]

A new extension of $\psi^{(m)}(q)$ has been
introduced in \cite{esmo2}, in connection with integrals involving
the polygamma and the loggamma functions.
These are the {\em balanced negapolygamma} functions, defined for $m \in
\allN$ by
\\
\ba\label{bal-negapolygamma}
\psi ^{( -m)} (q): = \frac{1}{m!}\left[{A_m (q)} - H_{m-1} B_m(q)\right].
\ea
\\
Here $H_{r}= 1 + 1/2 + \cdots + 1/r$ is
the harmonic number ($H_0=0$), $B_m(q)$ is the $m$-th Bernoulli
polynomial, and
\ba\label{Ak-def}
A_{m}(q) =  m\,\zeta'(1-m,q).
\ea
A function $f(q)$ is defined on $(0,1)$ is called {\em balanced} if its
integral over $(0,1)$ vanishes and $f(0) = f(1)$.
In \cite{esmo2} we have shown that
\ba\label{der-negapolygamma}
\frac{d}{dq}\psi^{(-m)}(q)=\psi^{(-m+1)}(q),\quad m\in\allN.
\ea

The function $\psi(z,q)$ defined in (\ref{psi-def1}) represents an extension
of these polygamma families to $q \in \allC$. Its main properties
are presented in the next theorem. The details appear in \cite{esmo4}. \\

\begin{Thm}\label{psi-main}
The generalized polygamma function $\psi(z,q)$ satisfies: \\

$\bullet$
For fixed $q \in \allC$, the function $\psi(z,q)$ is an entire
function of $z$.

$\bullet$ For $m \in \mathbb{Z} \; : \psi(m,q) = \psi^{(m)}(q)$.

$\bullet$ It satisfies
\ba
\frac{\partial}{\partial q} \psi(z,q) & = & \psi(z+1,q). \label{psi-der}
\ea

\end{Thm}

\section{Appendix 3: The relation between $Q_{j}$ and $R_{j}$}
\label{S:relations}

Using the relation (\ref{relation1}) we
can express the integrals $R_{j}$, defined by \eqref{def-R1}--\eqref{def-R6}
in terms of $Q_{j}$, defined by \eqref{def-Q1}--\eqref{def-Q6}. Recall that, for
$n> 1$ we have $h_{n} = 1+ 1/2 + \cdots + 1/(n-1)$ and $h_{1} = 0$.

{\allowdisplaybreaks
\begin{align}
\label{rel-RQ1}
R_{1}(n_{1},n_{2},n_{3}) & =  Q_{1}(n_{1},n_{2},n_{3}),\\
\label{rel-RQ2}
\pi R_{2}(n_{1},n_{2},n_{3}) & =  n_{3}! Q_{2}(n_{1},n_{2},n_{3}) +
h_{n_{3}} Q_{1}(n_{1},n_{2},n_{3}),\\
\label{rel-RQ3}
\pi^2 R_{3}(n_{1},n_{2},n_{3}) & =  n_{1}! n_{2}! Q_{3}(n_{1},n_{2},n_{3}) +
n_{2}! h_{n_{1}} Q_{2}(n_{1},n_{3},n_{2})\\
& +  n_{1}! h_{n_{2}} Q_{2}(n_{2},n_{3},n_{1}) +
h_{n_{1}} h_{n_{2}} Q_{1}(n_{1},n_{2},n_{3}), \nn \\
\label{rel-RQ4}
\pi^2 R_{4}(n_{1},n_{2},n_{3}) & =  n_{1}!n_{2}! Q_{4}(n_{1},n_{2},n_{3}) +
(-1)^{n_2} n_{1}! h_{n_{2}} Q_{2}(n_{2},n_{3},n_{1})\\
& + (-1)^{n_{1}+n_{3}} n_{2}! h_{n_{1}} Q_{2}(n_{1},n_{3},n_{2}) \nn\\
& + (-1)^{n_2} h_{n_{1}} h_{n_{2}} Q_{1}(n_{1},n_{2},n_{3}), \nn \\
\label{rel-RQ5}
\pi^3 R_{5}(n_{1},n_{2},n_{3}) & =  n_{1}!n_{2}!n_{3}! Q_{5}(n_{1},n_{2},n_{3}) +
n_{1}!n_{2}!h_{n_{3}}Q_{3}(n_{1},n_{2},n_{3}) \\
& + n_{1}!n_{3}!h_{n_{2}}Q_{3}(n_{1},n_{3},n_{2}) +
n_{2}!n_{3}! h_{n_{1}} Q_{3}(n_{2},n_{3},n_{1})  \nn \\
& + n_{1}!h_{n_{2}} h_{n_{3}}Q_{2}(n_{2},n_{3},n_{1}) +
n_{2}!h_{n_{1}}h_{n_{3}} Q_{2}(n_{1},n_{3},n_{2})  \nn \\
& + n_{3}!h_{n_{1}} h_{n_{2}} Q_{2}(n_{1},n_{2},n_{3}) +
h_{n_{1}} h_{n_{2}}h_{n_{3}} Q_{1}(n_{1},n_{2},n_{3}), \nn \\
\label{rel-RQ6}
\pi^3 R_{6}(n_{1},n_{2},n_{3}) & =  n_{1}!n_{2}!n_{3}! Q_{6}(n_{1},n_{2},n_{3}) +
(-1)^{n_3} n_{1}!n_{2}!h_{n_{3}}Q_{3}(n_{1},n_{2},n_{3}) \\
& + n_{1}!n_{3}!h_{n_{2}}Q_{4}(n_{1},n_{3},n_{2}) +
n_{2}!n_{3}! h_{n_{1}} Q_{4}(n_{2},n_{3},n_{1})  \nn \\
& + (-1)^{n_{3}} n_{1}!h_{n_{2}} h_{n_{3}}Q_{2}(n_{2},n_{3},n_{1})\nn\\
& + (-1)^{n_{3}} n_{2}!h_{n_{1}}h_{n_{3}} Q_{2}(n_{1},n_{3},n_{2})  \nn \\
& + (-1)^{n_{1}+n_{2}} n_{3}!h_{n_{1}} h_{n_{2}} Q_{2}(n_{1},n_{2},n_{3}) \nn\\
& + (-1)^{n_{3}} h_{n_{1}} h_{n_{2}}h_{n_{3}} Q_{1}(n_{1},n_{2},n_{3}) \nn .
\end{align}
}

\bigskip

\no
{\bf Acknowledgments}. The first author would like to thank the
Department of Mathematics at Tulane University for its hospitality, and
the partial support of grant MECESUP FSM0204. The
second author acknowledges the partial support of NSF-DMS 0070567.

\end{document}